\documentclass{ametsocV5}

\usepackage{amsmath,amsfonts,amssymb,bm}
\usepackage{mathptmx}
\usepackage{newtxtext}
\usepackage{newtxmath}
\usepackage{comment}
\usepackage{xcolor}
\usepackage{dirtytalk}

\title{Probability distributions for analog-to-target distances}

\authors{P. Platzer\correspondingauthor{Paul Platzer, paul.platzer@imt-atlantique.fr}}

\affiliation{Laboratoire des Sciences du Climat et de l'Environnement, Saclay, France \& Lab-STICC, UMR CNRS 6285, IMT Atlantique, F-29238, Plouzan\'{e}, France \& France \'{E}nergies Marines, Plouzan\'{e}, France}

\extraauthor{P. Yiou and P. Naveau}
\extraaffil{Laboratoire des Sciences du Climat et de l'Environnement (ESTIMR team), UMR 8212 CEA-CNRS-UVSQ, IPSL \& U Paris-Saclay, 91191 Gif-sur-Yvette, France}

\extraauthor{J-F. Filipot and M. Thiébaut}
\extraaffil{France \'{E}nergies Marines, Plouzan\'{e}, France}

\extraauthor{P. Tandeo}
\extraaffil{Lab-STICC, UMR CNRS 6285, IMT Atlantique, F-29238, Plouzan\'{e}, France}

\let\linenumbers\nolinenumbers\nolinenumbers
\let\internallinenumbers\nolinenumbers\nolinenumbers

\abstract{Some properties of chaotic dynamical systems can be probed through features of recurrences, also called analogs. In practice, analogs are nearest neighbours of the state of a system, taken from a large database called the catalog. Analogs have been used in many atmospheric applications including forecasts, downscaling, predictability estimation, and attribution of extreme events. The distances of the analogs to the target state condition the performances of analog applications. These distances can be viewed as random variables, and their probability distributions can be related to the catalog size and properties of the system at stake. A few studies have focused on the first moments of return time statistics for the best analog, fixing an objective of maximum distance from this analog to the target state. However, for practical use and to reduce estimation variance, applications usually require not just one, but many analogs. In this paper, we evaluate from a theoretical standpoint and with numerical experiments the probability distributions of the $K$-best analog-to-target distances. We show that dimensionality plays a role on the size of the catalog needed to find good analogs, and also on the relative means and variances of the $K$-best analogs. Our results are based on recently developed tools from dynamical systems theory. These findings are illustrated with numerical simulations of a well-known chaotic dynamical system and on 10m-wind reanalysis data in north-west France. A practical application of our derivations for the purpose of objective-based dimension reduction is shown using the same reanalysis data.}

\begin{document}


\maketitle

%








\section{Introduction}

Atmospheric analogs have been introduced by \cite{Lorenz1969} in a study on atmospheric predictability. The faster one target state $z$ and its best analog $a_0$ diverge from one another, the harder it is to predict the evolution of $z$. In Lorenz's study, the state $z$ was characterized by height values of the 200-, 500- and 850-mb isobaric surfaces at a grid of $\approx$ 1000 points over the Northern Hemisphere. The database of available analogs, called the catalog, contained five years of twice daily values. In his abstract, Lorenz states that there are ``numerous mediocre analogues but no truly good ones''.

Since Lorenz's work, analogs have been used in many applications such as weather generators \citep{Yiou2014}, data assimilation \citep{Hamilton2016, Lguensat2017}, kernel forecasting \citep{Alexander2017}, downscaling \citep{wetterhall2005statistical} climate reconstruction \citep{schenk2012reconstruction, fettweis2013important, yiou2013ensemble} and extreme event attribution \citep{cattiaux2010winter, jezequel2018role}.

The reason why Lorenz could not find any good analog was made clear later on by \cite{VanDenDool1994}. It was shown that for high-dimensional systems, the mean recurrence time of a good analog (identified as a minimum catalog size) grows exponentially with dimension. This result is a variant for analogs of the "curse of dimensionality", well known in data-sciences. With three pressure levels over the whole Northern Hemisphere, the dimension of Lorenz's study was very high, and only five years of twice-daily data was not enough to  hope  finding a good analog.

\cite{Nicolis1998} added a dynamical systems' perspective to Van den Dool's analysis. She showed that studying mean recurrence times was not enough, as the relative standard deviation of this recurrence time could be very high. Furthermore, it was shown that recurrence time statistics exhibit strong local variations in phase-space, so that certain target states may need a larger catalog size to find good analogs.

Accounting for Van den Dool's findings, it is now usual to reduce as much as possible the feature-space dimension before searching for analogs. Also, the last decades have witnessed a proliferation of data from in-situ and satellite observations, as well outputs from numerical physics-based model. Such conditions allows one to find good analogs in many situations, and it has become standard to use not just one, but many analogs (usually a few tens). From a statistical perspective, using many analogs instead of one can increase estimation bias, but it reduces estimation variance, so that the estimation is less sensitive to noise. Using many analogs also allows to perform local regression techniques on the analogs, such as local linear regression \citep{Lguensat2017}. This technique has proven efficient in analog forecasting applications \citep{Ayet2018}, and it was shown that local linear regression allows analog forecasting to capture the local Jacobian of the dynamics of the real system \citep{Platzer2020}.

This new context suggests to focus not only on the best analog $a_0$, but also the $k$-th best analog, for $k$ up to $\sim$40. Also, one can now reasonably hope to find good analogs using dimension reduction and a large amount of data.
Thus, one is less interested in recurrence times, but rather in analog performances. Performances of analog-based methods are largely conditioned by analog-to-target distances. In this work we propose to evaluate the probability distribution of these distances. Our analytical probability distributions make the link between analog-to-target distances, catalog size and local dimension. This brings new insight on the impact of dimensionality on analog forecasting performances. 

Section 2 outlines the theoretical framework and findings. The third section interprets the findings and compares this analysis with past studies. Section 4 shows results from numerical experiments of the \cite{Lorenz1963} system and from 10-m wind reanalysis data from the regional climate model AROME, further referred to as \say{the AROME reanalysis data}.

\section{Theory}\label{sec:theory}

\subsection{Analogs in dynamical systems and local dimensions}

For some dynamical system having an attractor set $\mathcal{A}$, (almost) all trajectories in the basin of attraction of $\mathcal{A}$ converge to the attractor \citep{milnor1985concept}. For such systems, almost all trajectories starting from the attractor come back infinitely close to their initial condition after a sufficiently long time \citep{poincare1890probleme}. Analog methods are based on the idea that if one is provided with a long enough trajectory of the system of interest, one will find analog states close to any point $z$ of the attractor $\mathcal{A}$.

The trajectory from which the analogs are taken is called the ``catalog'' $\mathcal{C}$, and can either come from numerical model output or reprocessed observational data. It can be seen either as a trajectory from a discrete dynamical system, or as evenly-spaced time samples from a continuous dynamical system. In any case, the catalog has a finite number of elements noted $L := \mathbf{card}(\mathcal{C})$. This catalog size may be divided by a typical correlation time-scale so that elements of the catalog can be considered independent \citep{VanDenDool1994}. In fact, for the analogs of a given target $z$ to be considered independent, it is enough that the typical distance between two analogs of $z$ be smaller than the typical distance between an analog and its time-successor.

The structure of the attractor, expressed by the system's invariant measure $\mu$, conditions the structure of the catalog and the ability to find analogs. In particular, \cite{VanDenDool1994} and \cite{Nicolis1998} studied the role of the attractor's dimension. Let $B_{z,r}$ the ball centered on $z\in\mathcal{A}$ and of radius $r$, then

\begin{equation}\label{eq:dzr}
    d_{z,r} := \frac{\log \mu (B_{z,r})}{\log r} \, ,
\end{equation}

\noindent defines the finite-resolution ($r$-resolution) local dimension at point $z$. Note that for ergodic measures, $\mu(B_{z,r})$ can be approximated by counting the number of times a given trajectory enters $B_{z,r}$ (this is the consequence of the ergodic theorem of \citealt{birkhoff1931proof}).

If $\mu$ is ergodic and $\lim_{r\to 0}d_{z,r}$ exists, then $\mu$ is said to be exact dimensional and the limit is independent of $z$ \citep{young1982dimension}. This typical value of the local dimension is noted $D_1$.

\begin{equation*}
    D_1 := \lim_{r\to 0}d_{z,r} \, .
\end{equation*}

The finite-resolution local-dimension $d_{z,r}$, however, can deviate from the typical value $D_1$. More precisely, $d_{z,r}$ exhibits large deviations from its limit value \citep[for more details, see][]{caby2019generalized}.

The distance from the $k$-th analog $a_k\in \mathcal{C}$ to the target state $z$ is noted $r_k := \mathrm{dist}(a_k,z)$. Distances are sorted so that $r_1 < r_2 < \cdots < r_K  $, and $K$ is the total number of analogs considered. Empirical methods usually set $K$ to a fixed value, reaching for a bias-variance trade-off. This amounts to looking at a lower quantile of the function $x\mapsto\mathrm{dist}(x,z)$. Another possibility is to set a threshold $R$ for the analog-to-target distances so that $r_K < R <r_{K+1}$. In this case, $K$ depends on $z$.

\subsection{Simple scaling of analog-to-target distance with dimension}\label{subsec:rkd_scaling}

Using extreme value theory and dynamical systems theory, \cite{caby2019generalized} showed that $d_{z,r}$ can be estimated using the empirical distribution of points inside a ball of exponentially decreasing radius. This empirical distribution is actually exactly the cumulative distribution function of the best available analogs. It then follows from \cite{caby2019generalized} that, for regular enough measures, we have the approximate scaling:

\begin{equation}\label{eq:rk_d}
	r_k(z) \sim k^{1/d},
\end{equation}

\noindent where $d=d_{z,r_K}$ is the local dimension at finite resolution $r_K$ (the largest analog-to-target distance). An application of this method to the three-variable system of \cite{Lorenz1963} is given in Fig. \ref{fig:loc_dim}.

Eq. \ref{eq:rk_d} gives an important point of our analysis, which is the scaling of $r_k$ with $k$, and is approximately given by a power-law with exponent $1/d$. However, this formula comes from a work on local dimensions, not analog-to-target distances. It is therefore not surprising that some of the elements required for our study are missing. In particular, this scaling does not give the constant in front of $k^{1/d}$, in which resides the relation to the catalog size, a crucial point for analog applications. Also, it only gives a mean or typical value of $r_k$, while our objective is to evaluate the probability distribution of $r_k$, or at least the probability of departures from this mean scaling.

The next section gives theoretical elements to evaluate the full probability distribution of $r_k(z)$ from the local dimension, the catalog size, and the analog number $k$.

\subsection{Full probability distribution of analog-to-target distance}

\subsubsection{Poisson distribution of the number of analogs in a ball}

\cite{Haydn2019} have shown that, for dynamical systems having Rare Event Perron-Frobenius Operator properties, and for non-periodic points $z$, the number of visits $k(z,r)$ of a trajectory of size $L$ into the ball $B_{z,r}$ follows a Poisson distribution with mean $L\mu(B_{z,r})$:

\begin{equation}\label{eq:Poisson}
	\mathbb{P}\left( k(z,r) = k \right) = \frac{\left(L\mu(B_{z,r})\right)^k}{k!}e^{-L\mu(B_{z,r})},
\end{equation}

\noindent where $k!$ is $k$ factorial. In the context of analogs, this is the probability to find $k$ analogs with distances to $z$ below the radius $r$. In the following we write $\mu_{z,r}:=\mu(B_{z,r})$.

\subsubsection{Distribution of analogs close to the sphere}

This section aims at using $\mu$ to evaluate $\mathbb{P}\left(r_k\in [r,r+\delta r)\right)$, the probability that the $k$-th analog-to-target distance is between $r$ and $r+\delta r$, for fixed $k$ and $z$ and where $\delta r$ is small compared to $r$.

The event ``$r_k\in [r,r+\delta r)$'' is the intersection of the event ``there are $k-1$ analogs in the ball $B_{z,r}$'' and the event ``there is one analog in $B_{z,r+\delta r}\cap \overline{B_{z,r}}$''. For a Poisson point process these two events are independent, so that:

\begin{equation}\label{eq:pr_rdr}
\begin{split}
\mathbb{P}\left( r_k\in [r,r+\delta r)\right) & = \mathbb{P}\left( \; k(z,r) = k-1  \; \land \;  \exists x \in \mathcal{C} \cap B_{z,r+\delta r} \cap \overline{B_{z,r}} \; \right) \\
& = \mathbb{P}\left( k(z,r) = k-1 \right) \; \mathbb{P}\left( \exists x \in \mathcal{C} \cap B_{z,r+\delta r} \cap \overline{B_{z,r}} \right) \\
& = \frac{\left(L\mu_{z,r}\right)^{k-1}}{(k-1)!}e^{-L\mu_{z,r}} \;  \mathbb{P}\left( \exists x \in \mathcal{C} \cap B_{z,r+\delta r} \cap \overline{B_{z,r}} \right) \, .
\end{split}
\end{equation}

Then, it follows from \cite{Haydn2019} that the event that strictly one element of the catalog lies between $B_{z,r}$ and $B_{z,r+\delta r}$ has a probability of the same form as Eq. (\ref{eq:Poisson}) but replacing $k$ by 1 and $\mu_{z,r}$ by $\delta \mu_{z,r} := \mu_{z,r+\delta r}-\mu_{z,r}$

\begin{equation}\label{eq:pr_!xdr}
\mathbb{P}\left( \exists ! x \in \mathcal{C} \cap B_{z,r+\delta r} \cap \overline{B_{z,r}} \right) = L \delta \mu_{z,r} e^{-L\delta \mu_{z,r}} \, .
\end{equation}

If the invariant measure $\mu$ is regular enough so that $\lim_{\delta r\rightarrow 0}\delta \mu_{z,r}=0$ we then have $e^{-L\delta \mu_{z,r}}\approx 1$. Also, the probability to find more than one element of the catalog between $B_{z,r}$ and $B_{z,r+\delta r}$ has a probability of $\mathcal{O}(\delta \mu_{z,r})^2$. This justifies the approximation $\mathbb{P}\left( \exists x \in \mathcal{C} \cap B_{z,r+\delta r} \cap \overline{B_{z,r}} \right) \approx \mathbb{P}\left( \exists ! x \in \mathcal{C} \cap B_{z,r+\delta r} \cap \overline{B_{z,r}} \right)$. Finally, combining Eq. (\ref{eq:pr_rdr}) and Eq. (\ref{eq:pr_!xdr}), one finds:

\begin{equation}\label{eq:p_rk_mu}
	\mathbb{P}\left( r_k \in [r,r+\delta r) \right) = L \delta\mu_{z,r} \frac{\left(L\mu_{z,r}\right)^{k-1}}{(k-1)!}e^{-L\mu_{z,r}} \; .
\end{equation}

This last equation is a more general form of our main result which is given in the next section. Here, the probability is expressed in terms of the invariant measure, which is usually not known analytically. The next section expresses the same probability in terms of the analog-to-target distance $r$.

\subsubsection{Distribution of analogs-to-target distance}

The link between $\mu_{z,r}$ and $r$ is given by the definition of the finite-resolution local dimension in Eq. (\ref{eq:dzr})

\begin{equation}\label{eq:mu_r}
\mu_{z,r} = r^d \, ,
\end{equation}

\noindent where $d=d_{z,r}$. The link between $\delta \mu_{z,r}$ and $\delta r$ involves variations of the local dimension wih $r$. Let $\Delta=d_{z,r+\delta r}-d_{z,r}$, we have:

\begin{equation}
\frac{\delta \mu_{z,r}}{\mu_{z,r}} = \left(1+\frac{\delta r}{r}\right)^{d+\Delta}e^{\Delta \log r} - 1 \, .
\end{equation}

Using the regularity hypothesis $\Delta\ll d$, and keeping only lower-order terms, we find:

\begin{equation}
\frac{\delta \mu_{z,r}}{\mu_{z,r}} = d\, \frac{\delta r}{r} + \Delta \log r \, .
\end{equation}

The term $d\, \frac{\delta r}{r}$ represents an almost steady increase in $\mu_{z,r}$ when $r$ grows. The term $\Delta \log r $ represents fluctuations in this increase given by the fluctuations in $d_{z,r}$. In practice, the method described in Sec. \ref{subsec:rkd_scaling} to evaluate $d$ should catch a mean local dimension over the analogs and not catch the fluctuations of $d_{z,r}$ with $r$ at scales smaller than $r_K$. Thus, the approximation :

\begin{equation}\label{eq:dmu_dr}
\frac{\delta \mu_{z,r}}{\mu_{z,r}} \approx d\, \frac{\delta r}{r} \, ,
\end{equation}

\noindent which is not valid in theory, should be relevant in practice for finite catalog size and regular enough measures. For small enough $\delta r$, one can then define $p_k$, the probability density function of $r_k$ through the identity $\mathbb{P}\left(r_k\in [r,r+\delta r)\right)=p_k(r)\delta r$. Combining Eq. (\ref{eq:p_rk_mu}), Eq. (\ref{eq:mu_r}) and Eq. (\ref{eq:dmu_dr}), we find:

\begin{equation}\label{eq:p_rk}
	p_k(r) = d\,L\,r^{d-1}\, \frac{\left(L\,r^{d}\right)^{k-1}}{(k-1)!} \, e^{-L\,r^{d}} \; .
\end{equation}

This last equation is our main result. An alternative proof for Eq. (\ref{eq:p_rk}) using extreme value theory is given in appendix A. Eq. (\ref{eq:p_rk}) then allows to compute the mean and variance of $r_k$ for fixed $k$ and $d$ :

\begin{subequations}\label{eq:rk_moments_exact}

	\begin{equation}
		\langle \, r_k \, \rangle \, = \, \frac{\Gamma\left(k+\frac{1}{d}			\right)}{L^{1/d}\; \Gamma(k)} \, ,
	\end{equation}

	\begin{equation}
		\langle \, r_k^2 \, \rangle \, - \, \langle \, r_k \, \rangle^2 \, = \, 				\frac{1}{L^{2/d}\; \Gamma(k)^2} \, \left\lbrace 			\Gamma\left( k+\frac{2}{d} \right)\Gamma(k) - 	\Gamma\left(k+\frac{1}{d}\right)^2 \right\rbrace \, ,
	\end{equation}

\end{subequations}

\noindent where $\Gamma$ is Euler's Gamma function. These identities can be simplified through scalings of the Gamma function $\Gamma(x+1) = \int_0^{+\infty} u^x e^{-u} \mathrm{d}u$ for large $x$, using Laplace's method up to second order to evaluate the integral (the first order gives Stirling's formula). This gives:

\begin{subequations}\label{eq:rk_moments_scaling}

	\begin{equation}
		k\geq 2 \, , \; \langle \, r_k \, \rangle \, \approx \, \left( \frac{k}{L} \right)^{1/d}\; ,
	\end{equation}

	\begin{equation}
		\frac{ \left(\langle \, r_k^2 \, \rangle \, - \, \langle \, r_k \, \rangle^2\right)^{1/2}}{\langle \, r_k \, \rangle} \, \approx \, \frac{1}{dk^{1/2}} \; ,
	\end{equation}

\end{subequations}

\noindent where we find again the scaling $r_k\sim k^{1/d}$ of Eq. (\ref{eq:rk_d}). These approximations will be increasingly valid as $k$ grows, but even for $k=2$, Eqs. (\ref{eq:rk_moments_scaling}a,b) give a satisfactory numerical approximation of Eqs. (\ref{eq:rk_moments_exact}a,b).

One can also compute $r^*_k$, the value of $r$ for which $p_k$ reaches a maximum:

\begin{equation*}
	r^*_k = \mathrm{argmax}_r\left\lbrace p_k(r) \right\rbrace = \left( \frac{k-\frac{1}{d}}{L} \right)^{1/d} \, ,
\end{equation*} 

\noindent and when $kd\leq1$, $r^*_k = 0$ and $p_k(0)=+\infty$. Note that the three quantities $\langle r_k\rangle$, $\left(\frac{k}{L}\right)^{1/d}$ and $r^*_k$ are equivalent as $k\rightarrow +\infty$.\medskip

\begin{figure}[t]
	\noindent\includegraphics[scale=.6]{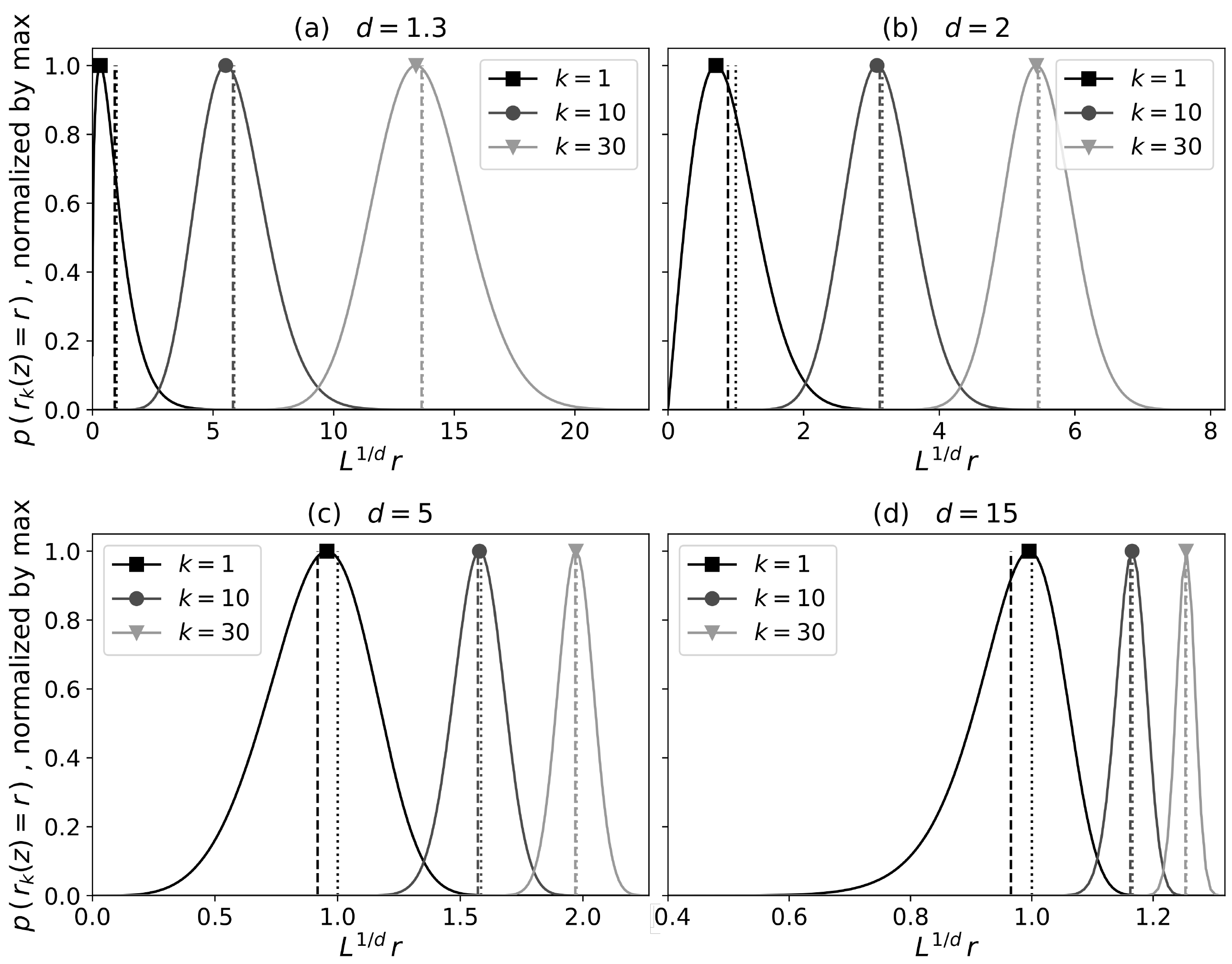}
	\caption{Probability density functions of $r_k$, the $k$-th analog-to-target distance, for fixed values of $k$, and of the local dimension $d$, from Eq. (\ref{eq:p_rk}). The dimension equals (a) 1.3, (b) 2, (c) 5, (d) 15. All densities $p_k$ are normalized by their maximum value. The distances are normalized by $L^{-1/d}$. Dashed vertical lines indicate the exact mean value $\langle r_k\rangle$ from Eq. (\ref{eq:rk_moments_exact}a), while dotted vertical lines indicate the approximate value $(k/L)^{1/d}$ from Eq. (\ref{eq:rk_moments_scaling}a). The argmax values of $p_1$, $p_{15}$ and $p_{30}$ are shown respectively with squares, circles and triangles.}
	\label{fig:proba_rk_r}
\end{figure}

Fig. \ref{fig:proba_rk_r} shows plots of $p_k(r)$ against $r$ for varying values of $d$ and $k$. As a consequence of the scaling $r_k\sim k^{1/d}$, we observe large variations of $\langle r_k\rangle$ with $k$ for small dimensions $d$, and very small variations of $\langle r_k\rangle$ with $k$ for large dimensions $d$. Note that, in the limiting case $d\to \infty$, the random variables $r_k$ are degenerate and all equal $L^{-1/d}$ almost surely. This can be witnessed through the different scales of the horizontal axis of the plots. Also, as a consequence of Eqs. (\ref{eq:rk_moments_scaling}), we have that the standard deviation of $r_k$ is a growing function of $k$ for $d<2$, while it is constant for $d=2$ and decreasing for $d>2$. However, the relative standard deviation of $r_k$ is always a decreasing function of $k$ and $d$ according to Eq. (\ref{eq:rk_moments_scaling}b).

\subsection{Rescaling and convergence to the standard Normal distribution}\label{subsec:rescale_normal}

Eqs. (\ref{eq:rk_moments_scaling}a,b) suggest the change of variables from $r$ to $u$ with

\begin{equation*}
    u = dk^\frac{1}{2}\left( \left(\frac{L}{k}\right)^{\frac{1}{d}}r-1 \right) \, ,
\end{equation*}

\noindent so that the probability density function of $u_k$, noted $h_k(u)$, is

\begin{equation}\label{eq:rescaled_true}
    h_k(u) = \frac{k^{k-\frac{1}{2}}}{(k-1)!} \left( 1 + \frac{u}{dk^\frac{1}{2}} \right)^{dk-1} \exp\left\lbrace -k\left( 1 + \frac{u}{dk^\frac{1}{2}} \right) \right\rbrace \, ,
\end{equation}

\noindent and simple asymptotic analysis gives

\begin{equation*}
    \lim_{k\rightarrow +\infty} h_k(u) = \frac{1}{\sqrt{2\pi}} \exp\left\lbrace-\frac{u^2}{2}\right\rbrace \, ,
\end{equation*}

\noindent which shows that the rescaled random variable $u_k$ converges in distribution to the standard Normal distribution as $k\rightarrow +\infty$. Note, however, that this limit should be hard to observe in practice, as the distribution of Eq. (\ref{eq:p_rk}) is valid only in the limit of large catalog size and with $k\ll L$.

\section{Consequences for applications of analogs}

\subsection{Comparison with previous studies}

The pioneering work of \cite{VanDenDool1994} focuses on the minimum length of catalog needed to have a 95\% chance to find at least one analog with a distance below a low threshold $\varepsilon$. With our notations, this condition can be written

\begin{equation*}
    L \; \mid \; \mathbb{P}(r_1<\varepsilon) > 0.95 \, .
\end{equation*}

\cite{VanDenDool1994} uses a Gaussian approximation for the difference between two states, which is reasonable in high dimensions. Then $\mathbb{P}(r_1<\varepsilon) = 1-(1-\alpha^{D_1})^L$, where $\alpha$ is the probability that the distance between two arbitrarily chosen states is less then $\varepsilon$ and can be expressed as the integral of a Gaussian probability density function. For small $\varepsilon$, $\alpha=\mathcal{O}(\varepsilon)$ and $\alpha^{D_1}\ll 1$. This finally suggests

\begin{equation}\label{eq:DoolExpr}
    L > \frac{\log 0.05}{\log(1-\alpha^{D_1})} \approx \frac{- \log 0.05}{\alpha^{D_1}} 
\end{equation}

Similar results can be found from Eq. (\ref{eq:p_rk}). Indeed, one has  $\mathbb{P}(r_1<\varepsilon) = \int_0^\varepsilon p_1(r)\mathrm{d}r = 1 - \left( \exp(-\varepsilon^d) \right)^L $, so that $\alpha \approx \varepsilon$. Here, $D_1$ is replaced by the local finite-resolution dimension $d$. Thus, our analysis encompasses the one of \cite{VanDenDool1994}.

\cite{Nicolis1998} extended the work of \cite{VanDenDool1994}. Interpreting Eq. (\ref{eq:DoolExpr}) in terms of mean return times and using the formula from \cite{kac1959probability}, she found an expression of mean return times using the identity $\mu_{z,r}\approx r^{D_1}$ and a mean velocity. This theoretical analysis includes neither variations in phase space of the return time, nor variability of the return time due to the variability of the catalog for fixed $L$. However, \cite{Nicolis1998} performed empirical estimates of such variations of the return time, shading light on the pitfalls of an analysis limited to mean return times.

In the present paper, the point of view switches from statistics of return times to statistics of analog-to-target distance, and is extended to the $K$ best analogs rather then just the first one. The full probability distribution of Eq. (\ref{eq:p_rk}) gives a detailed view of the variability of the process of searching for analogs.

\subsection{Searching for analogs: consequences}

The full probability distribution of Eq. (\ref{eq:p_rk}) has many consequences for the practical search of analogs.

For very low-dimensional systems ($D_1<2$), the first analog-to-target distance has a lower variability than the next ones, so that a given value of $r_1$ will be more representative of the next values of $r_1$ than a given value of $r_{10}$ would be of the next values of $r_{10}$. The inverse phenomenon happens for higher dimensional systems ($D_1>2$). 
This can be taken into account to evaluate the expected performances of analog methods.

Also, the scaling $r_k\sim k^{1/d}$ implies that the growth with $k$ of the mean analog-to-target distance is much faster for low-dimensional systems ($D_1\lesssim 2$), so that the 30-th analog would be much farther from $z$ than the first one. This would justify the use of a lower number of analogs $K$ in low-dimensional spaces, while high values of $K$ would not have a great impact on analog-to-target distances in high dimensions (see the abscissa of the lower-left panel in Fig. \ref{fig:proba_rk_r}).

For instance, \cite{Lguensat2017} use analogs to produce forecasts of several well-known dynamical systems, setting $K=40$, while the use of Gaussian kernels with a variable bandwidth equal to $\lambda_z = \mathrm{median}_k r_k$ allows to discard analogs with $r_k> \lambda_z$. One might think that the filtering out of analogs with $r_k> \lambda_z$ make the forecast procedure relatively insensitive to the choice of $K$. Conversely, assuming that $\lambda_z \approx\,  \langle r_{[K/2]} \rangle $ where $[K/2]$ is the integer part of $K/2$, we have that $\lambda_z$ grows with $K$ as $\lambda_z\sim K^{1/d}$. Thus, for low-dimensional systems such as the one of \cite{Lorenz1963} for which $D_1\approx 2.06$, our results suggest that high values of $K$ would have detrimental effects on the efficiency of analog methods.

Moreover, the scaling $\langle r_k \rangle \sim \left( \frac{k}{L} \right)^{\frac{1}{d}}$ can be used in the context of dimension reduction. Assume that one wants to perform a statistical task that necessitates $K$ analogs (for instance, an ensemble forecast). Then, assume that one wants to reduce the dimension in order to have $\langle r_K \rangle < \varepsilon$. From the scaling $\langle r_k \rangle \sim \left( \frac{k}{L} \right)^{\frac{1}{d}}$, we find that the dimension must be reduced to at least $d_{\mathrm{max}, K} = \left( 1 - \frac{\log(K)}{\log(L)} \right) d_{\mathrm{max}, 1}$. Detailed arguments and a practical example are given in Sec. \ref{sec:num}.\ref{subsec:dim_red}. Thus, for instance, if the criterion $\langle r_1 \rangle < \varepsilon$ is met for $d_{\mathrm{max}, 1}=10$ and if $L=10^4$, then the criterion  $\langle r_{25} \rangle < \varepsilon$ will be met only for $d_{\mathrm{max}, 25}=6$. This shows that any dimension reduction performed with the objective of increasing analog performances strongly depends on how many analogs are required.

Finally, the joint distribution of analog-to-target distances from appendix A theoretically allows to express the probability distributions of any random variable of the form $\sum_k \omega_k r_k^p $, where $\lbrace \omega_k \rbrace_k$ are weights and $p$ is a positive integer. Such quantities can give error bounds for analog methods (see \citealt{Platzer2020}, for the case of analog forecasting). However, a closed form for the distribution of such variables is yet to be derived.

\section{Numerical experiments}\label{sec:num}

\subsection{Three-variable Lorenz system}

\begin{figure}[t]
    \noindent\includegraphics[scale=.5]{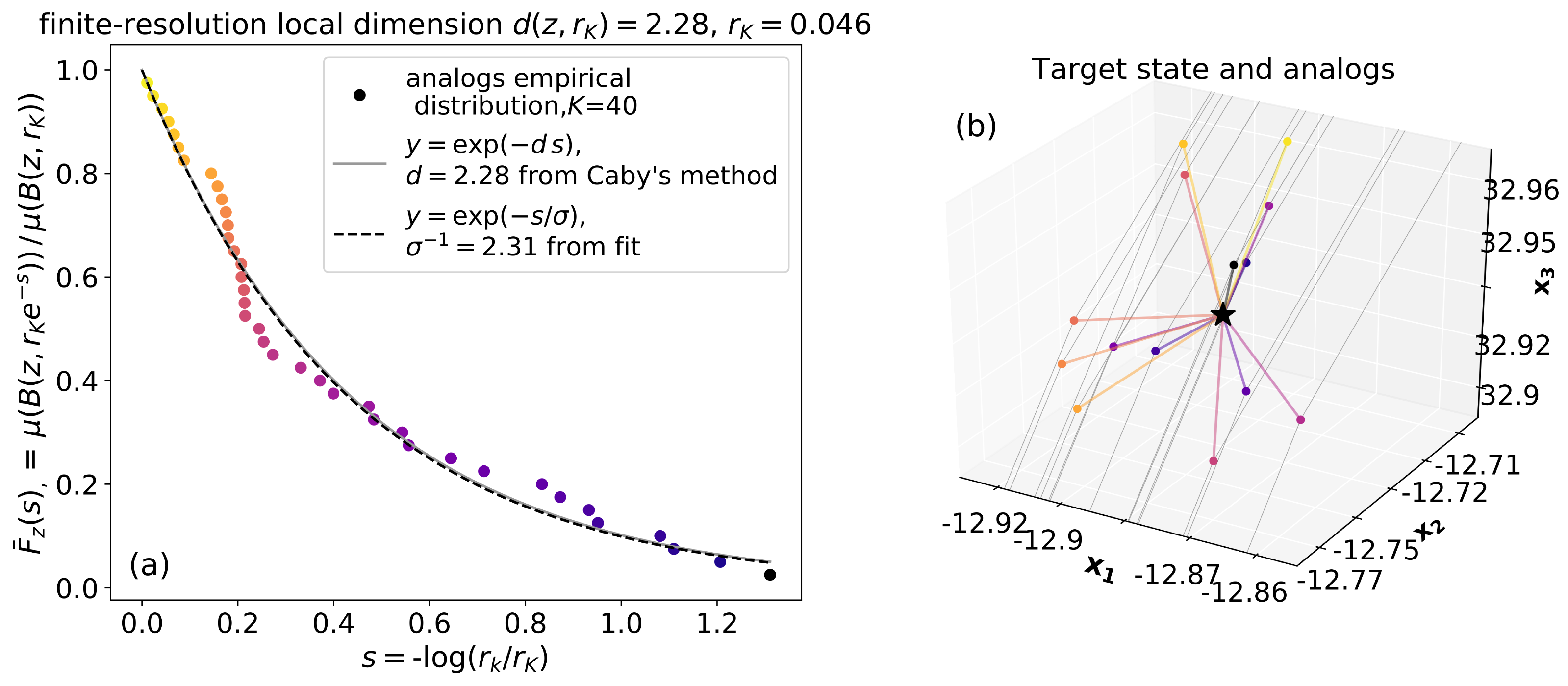}
    \caption{Computing the finite-resolution local dimension $d_{z,r_K}$ at a point $z$ of the three-variable \cite{Lorenz1963} system. (a) Following from \cite{caby2019generalized}, we evaluate $d$ by taking the mean of the empirical cumulative distribution function. For this example, fitting the empirical CDF with an exponential $\exp(-s/\sigma)$ and taking the inverse of $\sigma$ would have given approximately the same value for $d$. (b) Target $z$ (black star) and one in three analogs (colored dots matching (a)). The trajectories from which the analogs are taken are in grey. In this example, the typical analog-to-successor distance is much larger than the typical analog-to-target distance.}
    \label{fig:loc_dim}
\end{figure}

\begin{figure}[t]
	\noindent\includegraphics[scale=.45]{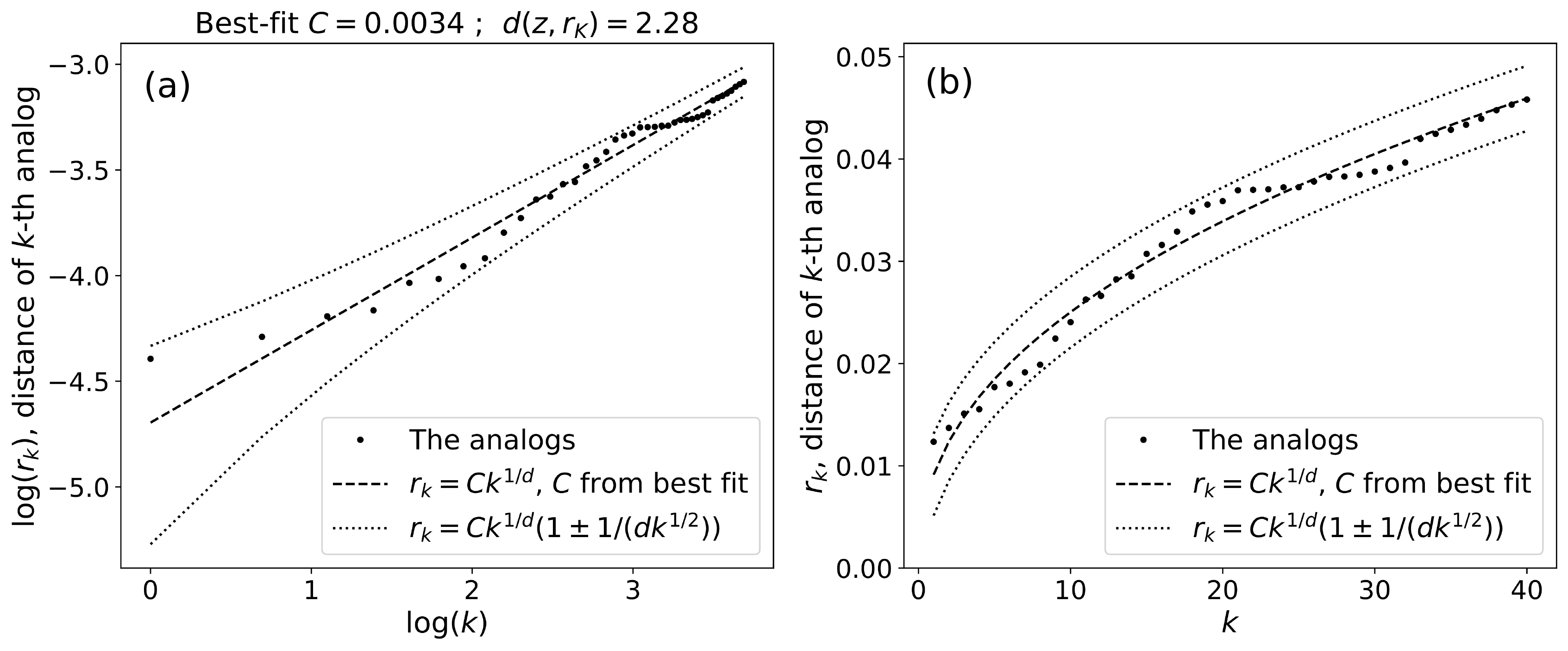}
	\caption{Analog-to-target distance $r_k$, against analog number $k$ at the same point $z$ than in Fig. \ref{fig:loc_dim}. (a) Log-scale. (b) Linear scale. Full circles are the empirical points given by the analogs. The dashed dark line is the best fit from equation (\ref{eq:rk_c}) where $d$ is fixed (from Caby's method) and $C$ is estimated with least-squares in log-scale. Assuming that this fit gives an estimation of the mean, the dotted lines represent approximate standard-deviation around this mean.}
	\label{fig:rk_power}
\end{figure}

Using the procedure of \cite{caby2019generalized}, one estimates the local finite-resolution dimension $d(z,r_K)$ for any point $z$ using the $K$-best analogs in the system of \cite{Lorenz1963}, hereafter noted L63. This procedure is illustrated in Fig. \ref{fig:loc_dim}. Then, the scaling of Eq. (\ref{eq:rk_moments_scaling}a) is used to make a least-squares fit from the data 

\begin{equation}\label{eq:rk_c}
 r_k(z) \approx^\mathrm{LS} C(z) k^{1/d}\, ,
\end{equation}

\noindent where $r_k(z)$ is the observed $k$-th analog-to-target distance and $\approx^\mathrm{LS}$ means that the constant $C(z)$ is evaluated with least-squares from Eq. (\ref{eq:rk_c}). Fig. \ref{fig:rk_power} shows an application of this procedure for a given $z$ of the L63, plotting the real values of $r_k$, and using $C(z) k^{1/d}$ as an approximation for $\langle r_k\rangle$ and dotted lines show the standard deviation around the mean from the approximate Eq. (\ref{eq:rk_moments_scaling}b).

From Eq. (\ref{eq:rk_c}) and Eq. (\ref{eq:rk_moments_scaling}) one expects to find:

\begin{equation}\label{eq:C_Ld}
	C(z) \approx L^{-1/d} \, ,
\end{equation}

\noindent however, as $L$ takes large values (from $10^5$ to $10^7$ or more), a small estimation error for $d$ results in a large estimation error for $L^{-1/d}$. Another way to look at this estimation issue is that $d$ is relatively insensitive to a rescaling of distances. Let:

\begin{equation}
    d' = \frac{\log \mu_{z,r}}{ \log (r/\rho)} \,
\end{equation}

\noindent where $\rho$ is a scalar value and $r/\rho$ is a rescaled version of $r$. Then $d'\sim d$ as long as $| \log \rho | \ll | \log r |$. In particular, the method of \cite{caby2019generalized} is insensitive to a rescaling, as it involves only ratios of distances (see the horizontal axis of Fig. \ref{fig:loc_dim}). Thus, Eq. (\ref{eq:C_Ld}) does not hold when $C$ and $d$ are determined as explained above. This is why $C(z)$ is rather evaluated through Eq. (\ref{eq:rk_c}), which allows one to find the rescaling:

\begin{equation}\label{eq:C_rho}
    C(z) = \frac{ \rho(z) }{ L^{1/d} } \, .
\end{equation}

Note that similar issues are raised by \cite{faranda2011extreme} regarding the continuity of $\mu_{z,r}$ with respect to $r$ and its limiting behaviour for small $r$, which motivates \cite{lucarini2014towards} to postulate that $\mu_{z,r}$ is the product of $r^{D_1}$ and a slowly varying function of $r$, which is in some sense equivalent to our hypothesis that $C(z)$ has to be rescaled when the local dimension is estimated from the method of \cite{caby2019generalized}.

\begin{figure}[t]
    \noindent\includegraphics[scale=.55]{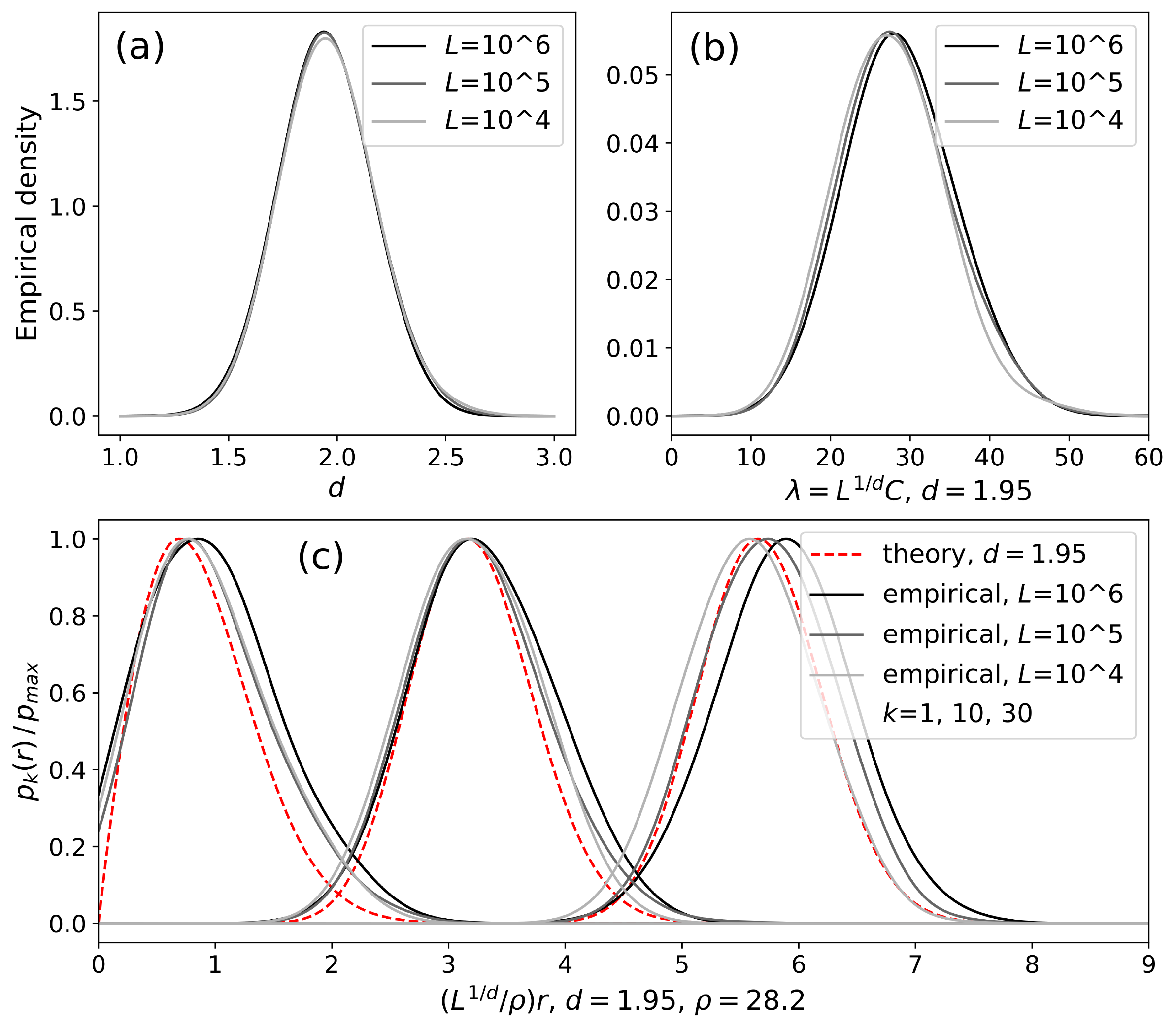}
    \caption{Numerical experiments of the system of \cite{Lorenz1963}, for a fixed target point $z$, using catalogs of various sizes $L$, repeating the experiment 600 times for each catalog to obtain empirical probability densities. (a) Empirical density of the local dimension $d$, obtained with the method of Fig \ref{fig:loc_dim} and with 150 analogs. (b) Empirical density of $\rho(z)$ obtained from Eq. (\ref{eq:rk_c}) and Eq. (\ref{eq:C_rho}), setting $d$ to the mean value of its empirical densities. (c) Normalized empirical probability densities of rescaled distances $\frac{L^{1/d}}{\rho}r$, setting $\rho$ and $d$ to the mean value of their empirical densities, and normalized theoretical probability densities using the same value of $d$. The probability densities are estimated using Gaussian kernels with bandwith of .15 (for $d$), 4 (for $\rho$), .3 (for rescaled $r$).}
    \label{fig:proba_rk_L63}
\end{figure}

Those formulas are tested in numerical experiments using the system of \cite{Lorenz1963}, with results reported in Fig. \ref{fig:proba_rk_L63}. Analogs of a fixed target point $z$ are sought for in $3\times 600$ independent catalogs, with three different catalog sizes. Each catalog is built from a random draw without replacement of $L$ points inside a (common) trajectory of $10^9$ points, generated using a Runge-Kutta numerical scheme with a time step of 0.01 in usual non-dimensional notations. The dimension is calculated using $K=150$ points, where this number is justified by a bias-variance trade-off: using this number and testing the procedure on 100 points picked from the measure $\mu$, one finds a mean dimension between 2.03 and 2.04, which is coherent with values reported by \cite{caby2019generalized}, and a standard deviation of $\sim0.26$. Using a lower value of $K$ results in a higher variance, and using higher values results in biases that are dependent on the value of $L$ used in this study.

The consistency of empirical densities of $\rho$ across varying values of $L$ validates the scaling of $C$ with $L$ and $d$. Empirical probability densities of rescaled analog-to-target distances, also consistent across varying catalog sizes, are coherent with the theoretical probability densities from Eq. (\ref{eq:p_rk}). The values of the rescaling parameter $\rho$ are not surprising, as typical values of distances between points in the attractor are $\sim 16$ and maximum distances are $\sim 28$. Note that \cite{Nicolis1998} uses a rescaling in studying analog return times with Lorenz's three-variable system, dividing all distances by the maximum distance between two points on the attractor.

Repeating this experiment for different target points $z$ gives similar results. Values of $\rho$ are of the same order of magnitude as the one reported in Fig. (\ref{fig:proba_rk_L63}). The consistency across varying values of $L$ is almost always recovered, except for some points that have slightly higher dimensions $d\gtrsim 2.15$ (not shown here). We expect this to come from a bad choice of $K$ when estimating the dimension and the rescaling factor: the choice of $K=150$ is relevant for most points, but should be adapted to the local dimension.

\subsection{AROME reanalysis data: dimensionality}

To further appreciate the applicability of our results to high-dimensional, real geophysical systems, the theoretical developments from Sec. \ref{sec:theory} are tested on five years (2015-2019) of hourly 10m-wind output from the physical model AROME \citep{ducrocq2005projet} coupled with satellite, radar, and in-situ observations through a variational data assimilation scheme \citep[similar to the one of][]{fischer2005overview}. The spatial domain is an evenly spaced grid above Britanny, with latitudes ranging from 47.075$^\circ$ to 49.3$^\circ$ and longitudes from -5.7$^\circ$ to -2.575$^\circ$, and a spacing of 0.025$^\circ$. To focus on wind at sea, land points are removed from the data resulting in a domain of 8190 grid points: this last step allows for comparison with ongoing work targeted at offshore wind characterization and forecast.

From this data, one can compute local dimensions with the method of \cite{caby2019generalized}. As the data is limited ($\sim 3 \times 10^4$ time-points), $K$ is set to 40. Note also that, as elements of the catalog are only one hour away from each other, they cannot be assumed independent. Therefore, if several analogs are neighbours in time, only one analog is retained, and it is selected randomly in the set of time-neighbouring analogs. Also, analogs that are less than one and a half days away from the target state $z$ are discarded.

\begin{figure}[t]
	\noindent\includegraphics[scale=.57]{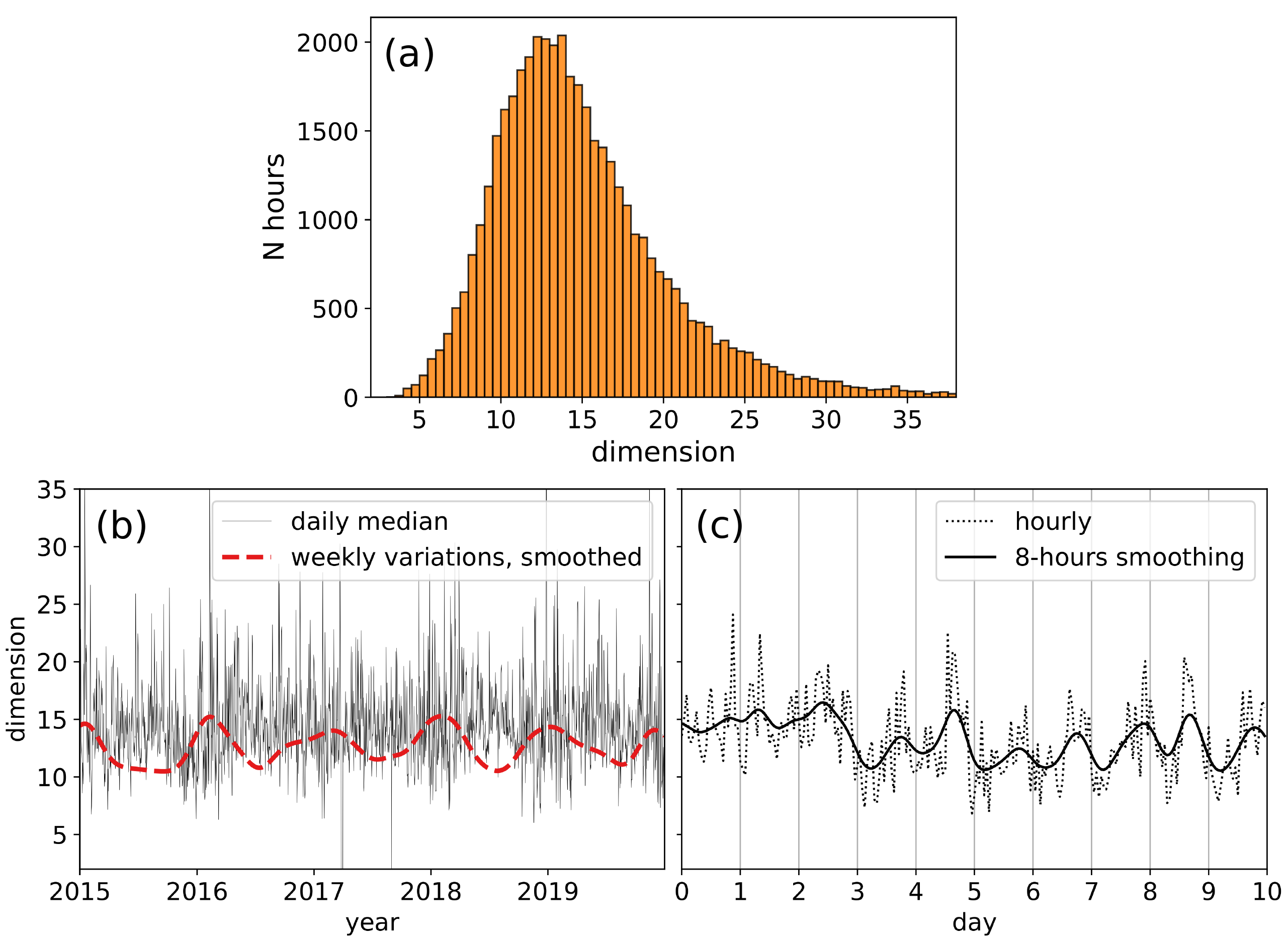}
	\caption{Statistics of local dimensions estimated using the method of \cite{caby2019generalized}, as in Fig. \ref{fig:loc_dim}. (a) Histogram of dimension from 10m-wind data off the Britanny coast. (b) Five years of dimension daily averages, and weekly variations defined as the difference between the 90\% and 10\% quantiles of hourly dimension over a week. This last quantity is smoothed over a $\sim$80-days window using convolution and Gaussian kernels. (c) Fourteen days of hourly local dimension, and a 8-days smoothing using convolution and Gaussian kernels.}
	\label{fig:d_winds}
\end{figure}

Histograms of local dimensions are plotted in Fig \ref{fig:d_winds}(a). These indicate that the system lives in an attractor of dimension approximately between 7 and 19, with some local dimensions likely to exceed 25. Our local dimension histogram is similar in shape to the one of \cite{faranda2017dynamical}, who also focused on North-Atlantic circulation. However, our histogram shows slightly higher average dimensions and a higher variability. Note that we focus on two components of horizontal wind velocity, on a dense grid of $\sim 10^4$ grid points, while \cite{faranda2017dynamical} focus on sea-level pressure (SLP) at $\sim 10^3$ grid points. Therefore, it is not surprising that we find higher average values of the local dimension. The fact that we observe a higher variability in the local dimension could be due to an intrinsic higher variability of this dynamical indicator, but also to a higher variability in the process of estimating $d$ caused by a lack of data. Indeed, we have slightly less data than \cite{faranda2017dynamical}, for a system of slightly higher dimension, so that we can find less good analogs to estimate $d$ than \cite{faranda2017dynamical}. \cite{faranda2017dynamical} use $L\sim 2 \cdot 10^4$ days of historical data. We use $\sim 4 \cdot 10^4$ hours of data, which must be divided by the typical correlation time-scale in hours. If we assume that the latter is between 12 and 24 hours, we find that our $L$ is between $1.5\cdot 10^3 $ and $3\cdot 10^3$. 

\cite{faranda2017dynamical} found a seasonality in the local dimension of SLP fields, with higher dimensions and a higher variability in winter. In our case, no seasonal trend for the mean or median dimension is observed, but the temporal variability of local dimensions is higher in winter, as witnessed in Fig. \ref{fig:d_winds}(b). Also, a diurnal cycle can be seen in Fig. \ref{fig:d_winds}(c), with dimension increasing in daytime and decreasing in nighttime. As diurnal variability is mixed with other sources of variability, it cannot always be identified by eye (see the 3 first days of Fig. \ref{fig:d_winds}c). Histograms of dimension restricted to daytime are similar to histograms restricted to nighttime, so that diurnal cycle does not appear to be the main driver of dimension variability.




This study of the dimension of the AROME reanalysis data using analogs can be compared with another method designed to categorize 10m-wind fields into classes, hereafter referred to as clusters. They are composed of the model grid points forming a geographic area of several thousands of square kilometers. We first adopted the empirical orthogonal functions (EOFs) approach to resolve separately the different spatial modes of the zonal and meridional wind velocities. We restricted the study to the 50 first EOFs which explain 98.9$\%$ and 98.7$\%$ of the total variance of the zonal and meridional velocities respectively. These EOFs were employed so as to compose a dataset of 8190 $\times$ 100 values (50 EOFs, 2 velocity components) used to feed a Gaussian mixture model (GMM, see \citealt{reynolds2009gaussian}).

GMM requires to impose as input the number of components, i.e. the number of clusters, in the model.  The optimum number of clusters was determined through the calculation of the Bayesian information criterion (BIC) score \citep[see][]{gideon1978estimating}. This score allows to select an optimal model to fit a dataset with a reasonable number of components. Low BIC scores correspond to a trade-off between the ability of the model to predict the data (the likelihood), and the number of model components. The relationship between the optimum number of clusters and the numbers $N$ of EOFs selected to feed the GMM was investigated. $N$ was set to be ranging from 10 to 100 with a step of 5. Our results (not shown) show that the optimum number of clusters is a slowly decreasing function of $N$. However, a convergence towards an optimum numbers of clusters of 10 is found when using the 50 (and more) first EOFs of the zonal and meridional velocities. Moreover, lower values of $N$ ($\leq$ 35) exhibited a less pronounced inflexion point of the BIC score curve, thus yielding a higher variability of estimation of the optimal number of clusters.

\begin{figure}[t]
    \noindent\includegraphics[scale=.45]{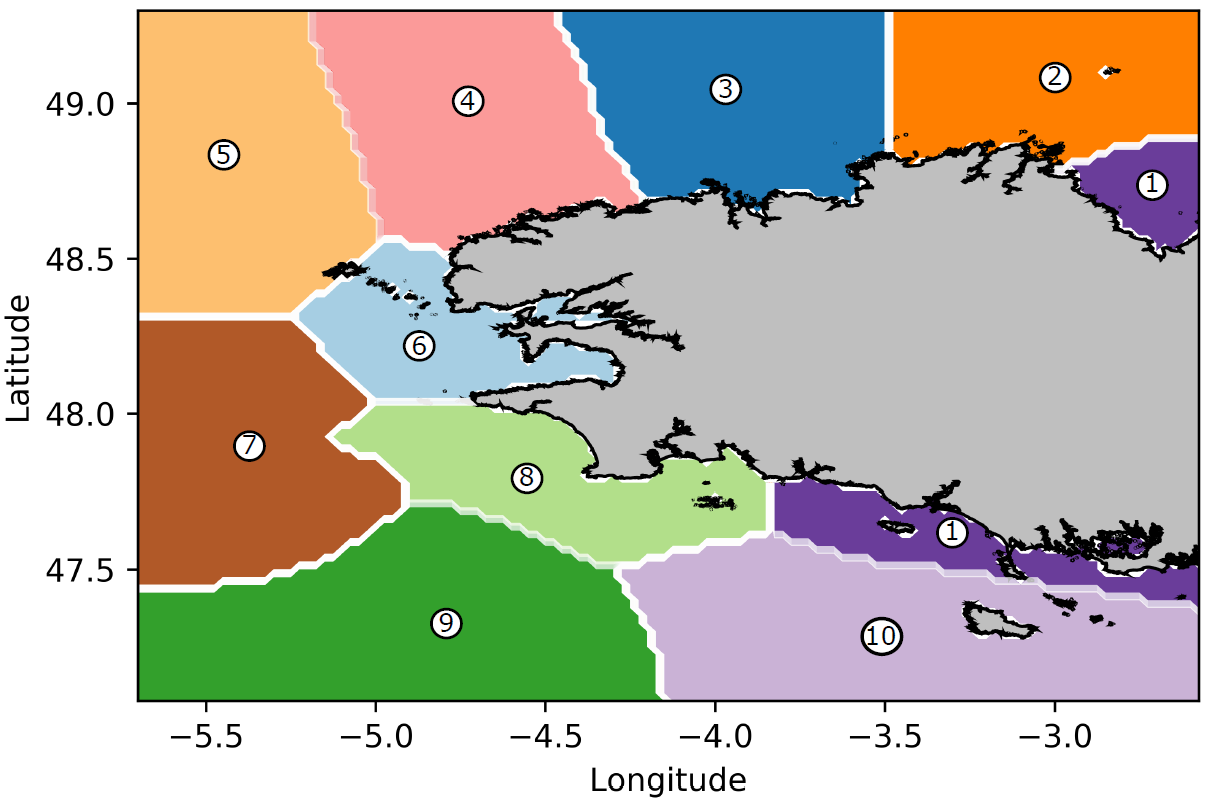}
    \caption{Spatial distributions of the wind clusters (colored area) blowing off the Brittany coast. The clustering was done through the use of a GMM. For ease of reading, the clusters are numbered from 1 to 10. The dataset employed for performing the GMM was composed of the 50 firsts empirical orthogonal functions (EOFs) of the zonal and meridional velocities derived from five years (2015-2019) of hourly 10m-wind output from the physical model AROME.}
    \label{fig:clusters}
\end{figure}

Fig. \ref{fig:clusters} shows the spatial distributions of the ten clusters. Cluster number 1 covers small areas on both sides of Brittany, at the eastern extremity. North of Brittany, the surface areas covered by the clusters (numbered from 2 to 5) are globally similar. Off the western extremity of the Brittany coast, the wind dynamics are more complex since the clusters 6 and 8 are associated with the smallest areas. South of Brittany, two clusters (9 and 10) are sufficient to explain the wind dynamics. Globally, the clusters connected to the coast are less spread out than the clusters located off the coast. This reveals the complexity of the land-sea transition which is associated with complex thermodynamic processes. 

The optimum number of clusters of 10 is in agreement with the dimensionality study using analogs, which indicate that the average attractor dimension is close to 13. This cluster study complements and strengthens our dimensionality analysis.

\subsection{AROME reanalysis data: analog distances}

\begin{figure}[t]
	\noindent\includegraphics[scale=.7]{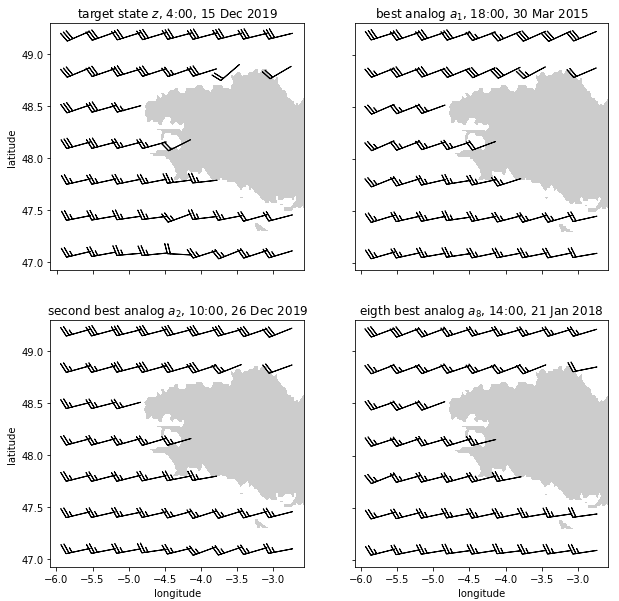}
	\caption{An example of target state $z$ and first, second and eight analogs, using 10m-wind data off the coast of Britanny from the AROME reanalysis. Standard station model notations are used, with wind speed in knots and point-centered flags.}
	\label{fig:ex_wind_ana}
\end{figure}

An example of target state and analogs is shown in Fig. \ref{fig:ex_wind_ana}. The chosen target state is a classical winter situation in Brittany, with strong eastward wind coming from the sea. Thus, good analogs are found in the catalog. It is hard to discriminate which analog is best: for such a high-dimensional system, the first analog-to-target distances are very similar.

By rescaling the variable $r_k$ in the following way:

\begin{equation}\label{eq:rescale_C}
	r_k \, \rightarrow \, d k^{1/2} \left( \frac{r_k}{Ck^{1/d}}-1 \right)\, , 
\end{equation}

\noindent where $d(z)$ is determined through the method of \cite{caby2019generalized} and $C(z)$ through the least-squares approximation introduced in the previous section, one should find probability distributions approaching a standard Normal distribution, as shown in Sec. \ref{sec:theory}.\ref{subsec:rescale_normal}, especially for large values of $k$. However, due to the small catalog size, only probability densities up to $k=8$ will be studied, otherwise the expressions obtained theoretically in the limit $L\to \infty$ are likely not to hold.

\begin{figure}[t]
    \noindent\includegraphics[scale=.58]{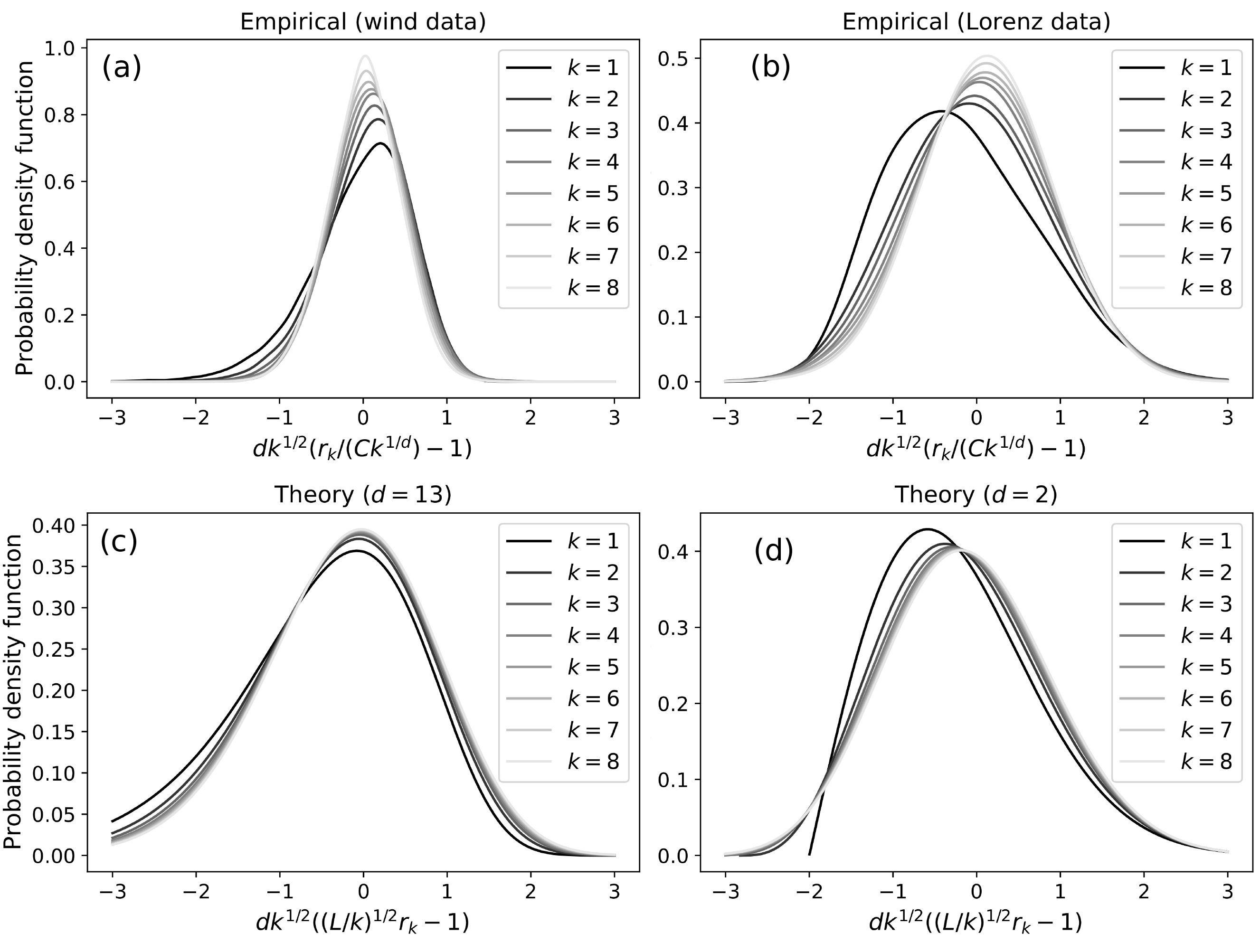}
    \caption{Probability densities of rescaled analog-to-target distances $r_k$ from 10m-wind data off the Britanny coast (a), and from numerical experiments of the Lorenz (\citeyear{Lorenz1963}) system (b), compared to theoretical distributions from Eq. (\ref{eq:rescaled_true}) for a local dimension of 13 (c) and 2 (d). Empirical probability densities are estimated using Gaussian kernels with a bandwidth of 0.3.}
    \label{fig:rk_winds}
\end{figure}

To obtain these distributions, analogs of each hourly $z\in \mathcal{C}$ (where $\mathcal{C}$ is the catalog) are sought for in the catalog, omitting analogs that are neighbours in time as explained previously. For each $z$, $C(z)$ is computed from Eq. (\ref{eq:rk_c}), and the distances are rescaled following Eq. (\ref{eq:rescale_C}) and then stored. Finally, the stored values of each rescaled $r_k$ are used to estimate probability density functions using Gaussian kernels with a bandwidth of 0.3. Fig. \ref{fig:rk_winds} shows the outcome of this procedure. For comparison, a similar procedure is applied on data from the model of \cite{Lorenz1963}, using a catalog of $L=10^6$ points and testing the procedure on $10^5$ target points that are taken from a trajectory independent from the catalog. Also, the theoretical density functions from Eq. (\ref{eq:rescaled_true}) are shown for similar (fixed) dimensions.

Fig. \ref{fig:rk_winds} shows a relatively good agreement between theoretical and empirical distributions, especially for the Lorenz data. Indeed, the curves of panels (b) and (d) are similar in shape, especially the asymmetry for $k=1$. As $k$ grows, the variance of the empirical data (b) becomes smaller than expected in theory (d). This can be explained by the fact that the assumption $L\to+\infty$ (or equivalently $r_k\to0$) is better satisfied for low values of $k$. High values of $r_k$ are associated with a low variability. This also explains the lower variance of the empirical curves (a) compared to the theoretical curves (c), using the wind data. Again, the asymmetry in the shape of the curves for $k=1$ is respected, and the estimation of the mean fits our theory.

\subsection{AROME reanalysis data: objective-based dimension reduction}\label{subsec:dim_red}

In this section, we assume that we want to reduce the dimension of the AROME reanalysis data in order to achieve the following criterion:

\begin{equation}\label{eq:criter}
    \frac{\overline{r_k}}{\mathrm{RMSD}} < \varepsilon \, 
\end{equation}

\noindent where $\overline{r_k}$ is the mean over all target points of the $k$-th analog-to-target distance, RMSD is the root-mean-squared distance between two points randomly taken from the dataset, and $\varepsilon$ is a user-defined threshold.  $\overline{r_k}$ is thus different from $\langle r_k(z) \rangle$, which is the mean over all possible realizations of the catalog, at fixed target $z$.

We assume that dimension is reduced using EOFs, which allows to reduce $\frac{\overline{r_k}}{\mathrm{RMSD}}$. However, one might not want to reduce dimension too much, in order to keep enough information on the state of the system. In this scenario, the practical question is: what is the maximum number of EOFs that can be used in order to meet Eq. (\ref{eq:criter}) ?

Following from the theoretical results of this paper, we assume that, for each target point $z$,

\begin{equation*}
    \frac{\langle r_k(z) \rangle}{\mathrm{RMSD}}  = \rho(z) \left( \frac{k}{L} \right)^{\frac{1}{d(z)}}\, ,
\end{equation*}

\noindent where $\rho(z)$ is of the order of 1. When using the method described in the previous sections to compute $d(z)$ and $C(z)$, we find that $\rho(z)$ is typically between 0.4 and 0.7. Then, we make the following ergodicity hypothesis:

\begin{equation*}
    \overline{r_k(z)} =  \overline{\langle r_k(z)\rangle} \, 
\end{equation*}

\noindent adding the hypothesis that $d(z)$ is almost constant, we finally find:

\begin{equation*}
    \frac{\overline{r_k}}{\mathrm{RMSD}} \approx \overline{\rho} \left(\frac{k}{L} \right)^{\frac{1}{d}} \, ,
\end{equation*}

\noindent which gives, combined with Eq. (\ref{eq:criter}):

\begin{equation*}
    d < d_{\mathrm{max},k} := \frac{\log(L/k)}{-\log(\varepsilon / \overline{\rho})} \, .
\end{equation*}

From this formula, it appears that $d_{\mathrm{max},k}$ is a linear function of $\log(k)$. This can be rearranged to give:

\begin{equation*}
    d_{\mathrm{max},k} = d_{\mathrm{max},1} \left( 1 - \frac{\log(k)}{\log(L)} \right) \, .
\end{equation*}

This last expression shows how $d_{\mathrm{max},k}$ strongly depends on $k$. On a practical example, assume that $d_{\mathrm{max}, 1}\approx 10$ and that $L=10^4$, then $d_{\mathrm{max}, 25}\approx 6$. Reducing dimension in order to improve the performance of analog methods thus strongly depends on how many analogs are needed for the analog method.

\begin{figure}[t]
    \noindent\includegraphics[scale=.7]{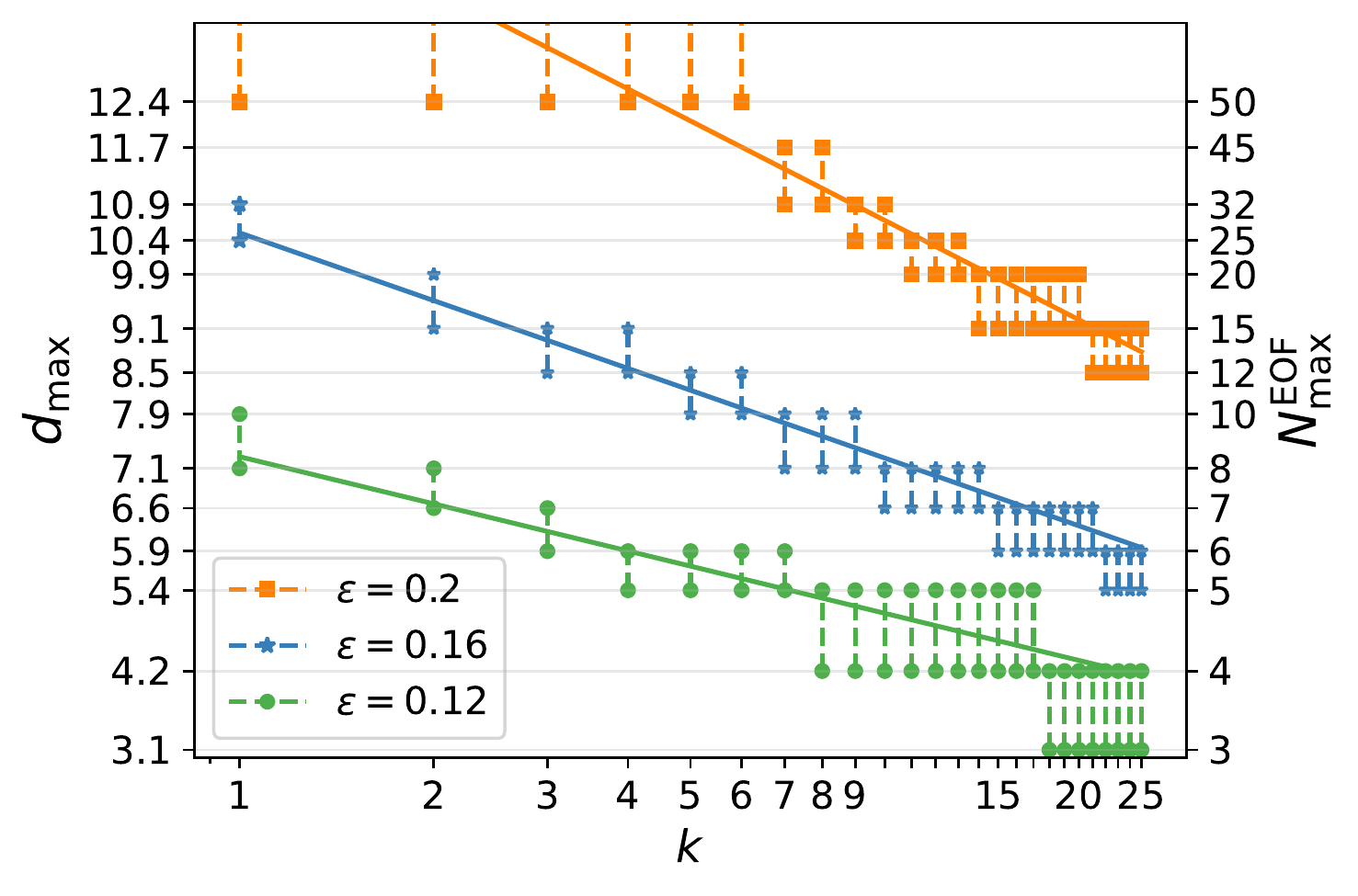}
    \caption{Maximum dimension (or number of EOF) to fulfill the criterion $ \frac{1}{\mathrm{RMSD}} \overline{r_k}  < \varepsilon  $, where $\overline{r_k}$ is the mean over all target points of the $k$-th analog-to-target distance, RMSD is the root-mean-squared distance between two random points from the data set, and $\varepsilon$ is a user-defined threshold. We use the 10m-wind data, and we project both component simultaneously on $N^\mathrm{EOF}$ basis functions. For a given value of $N^\mathrm{EOF}$, the dimension is computed as the mean of dimensions estimated from the method of \cite{caby2019generalized}. Then $ \frac{1}{\mathrm{RMSD}} \overline{r_k}$ is computed empirically, giving upper- and lower-bounds for the maximum dimension (or number of EOF). Full lines show the theoretical scaling $d_{\mathrm{max},k}=d_{\mathrm{max},1}\left( 1- \log(k)/\log(L) \right)$. The values of $d_{\mathrm{max},1}$ were set by hand on a visual criterion, and $L$ was set to $\sim 2\cdot 10^3$ which corresponds to a correlation time-scale of 24 hours.}
    \label{fig:dmax_k}
\end{figure}

Fig. \ref{fig:dmax_k} shows comparison of this scaling with numerical experiments performed on the AROME reanalysis data. Given the number of approximations that we have taken, the agreement between our theoretical scaling and the numerical experiments is satisfying.

\section{Conclusion}

We combined extreme value theory and dynamical systems theory to derive analytical joint probability distributions of analog-to-target distances in the limit of large catalogs. Those distributions shed new light on the influence of dimension in practical use of analog. In particular, we found that the number of analogs used in empirical methods induces larger biases in low dimension than in high dimension. Contrarily to previous works on the probability to find good analogs, this study focuses on distances rather than return times, and gives whole probability distributions rather than first moments. Numerical simulations of the three-variable Lorenz system confirm the theoretical findings. 10m-wind reanalysis data from the AROME physical model show that our analysis is also relevant for real systems. Our investigation indicates that the studied wind fields lie in an attractor of moderately high dimension $\sim$13, which is in agreement with the optimal number of regional cluster found using a Gaussian mixture model and the Bayesian information criterion. In this situation of moderate dimensionality, the analog-to-target distances of the first analogs are all very similar and have a low variability. Our theoretical derivations can be used to find optimal dimension reduction for improving analog performances, which we demonstrate on an example using the AROME reanalysis data.



\acknowledgments
The work was financially supported by ERC grant No. 338965-A2C2 and ANR No. 10-IEED-0006-26 (CARAVELE project). This piece of work took its origins in discussion with Théophile Caby, to whom we express our gratitude. The theoretical derivations of the probability density functions shown in this paper are the result of several exchanges with Benoît Saussol, who we must thank here. We are indebted to Fabrice Collard, Bertrand Chapron, and Caio Stringari, for fruitful insights and discussions about the exploration and interpretation of the AROME reanalysis data.

%
%


%






%
%
%

\appendix[A]\label{app:alt_proof}

\appendixtitle{Alternative proof for $p_k(r)$ and joint probability distribution}

\cite{Lucarini2016} give a detailed analysis of the map from $\mathcal{A}$ to $\mathbb{R}$, $x\mapsto-\log\mathrm{dist}(z,x)$, using tools from dynamical systems theory and extreme value theory (EVT). For our purpose, it is interesting to look at the simpler distance map $x\mapsto \mathrm{dist}(z,x)$.

The minimum of this map over the catalog is achieved for the best analog of $z$, $a_1$. The minimum is thus $r_1$. EVT tells \citep[see][]{coles2001introduction} that in the limit of large catalog, the minimum of this lower-bounded distance map on a finite sample of the attractor (a catalog of size $L$) follows a Weibull distribution, after rescaling. The Poisson law from Eq. (\ref{eq:Poisson}) with $k=1$ actually gives the scaling and the exact form of the Weibull distribution:

\begin{equation*}
\mathbb{P}(r_1>r)=e^{-Lr^d} \, ,
\end{equation*}

\noindent for positive $r$, otherwise the probability is 1.

The $K$ largest order statistics of this function then correspond to the $K$ analogs of the point $z$. Again, in the limit of large catalog and for small enough $K$, EVT provides the limit law \citep[see][]{coles2001introduction} for the $k$-th minima of this distance function when $L\to \infty$ :

\begin{equation*}
\mathbb{P}(r_k>r) = e^{-Lr^d} \sum_{s=0}^{k-1}\frac{(Lr^d)^s}{s!} \, .
\end{equation*}

Differentiating and with a bit of rearrangement, one finds back the formula of Eq. (\ref{eq:p_rk}).

\begin{equation*}
\begin{split}
p_k(r) & = -\frac{\partial}{\partial r}\mathbb{P}(r_k>r) \\
& = d\,L\,r^{d-1}\, \frac{\left(L\,r^{d}\right)^{k-1}}{(k-1)!} \, e^{-L\,r^{d}} \, .
\end{split}
\end{equation*}

From a broader perspective, extremal process theory \citep{lamperti1964extreme} gives the joint distribution of analog-to-target distances $p_{1:K}$ in the limit $L\to \infty$:

\begin{equation*}
	p_{1:K}(r_1,\ldots,r_K) = (dL)^K \left(\prod_{k=1}^K r_k \right)^{d-1} e^{-Lr_K^d} 
\end{equation*}

\noindent where the function is non-zero only when $0<r_1<r_2<\ldots<r_K$. For notation convenience and only in this formula, the random variables $r_k$ are noted identically as the values they can possibly take.

\appendix[B]

\appendixtitle{Three-variable Lorenz system}\label{app:lorenz}

The three-variable "L63"  \cite{Lorenz1963} system of equations is:

\begin{equation}\label{eq:L63}
\begin{cases}
 \dfrac{\mathrm{d}x_1}{\mathrm{d}t}=\sigma(x_2-x_1) \, ,\\ \\
 \dfrac{\mathrm{d}x_2}{\mathrm{d}t}=x_1(\rho-x_3)-x_2 \, ,\\ \\
 \dfrac{\mathrm{d}x_3}{\mathrm{d}t}=x_1x_2-\beta x_3 \, ,
 \end{cases}
\end{equation}

\noindent with usual parameters $\sigma=10$, $\beta=8/3$ and $\rho=28$.


\bibliographystyle{ametsoc2014}
\bibliography{biblio}

%

\end{document}